\documentclass{amsart}
\usepackage[dvips,final]{epsfig}
\usepackage{amsmath}
\usepackage{amscd}
\usepackage{amsfonts}
\usepackage{amssymb}
\usepackage{latexsym}
\title{
Around distance-squared mappings
}
\author{Shunsuke Ichiki}
\address{Dai Nippon Printing Co., Ltd.,
Tokyo 162-8001, Japan}
\email{ichiki-shunsuke-jb@ynu.jp}
\author{Takashi Nishimura
}
\address{Research Group of Mathematical Sciences,
Research Institute of Environment and Information Sciences,
Yokohama National University,
Yokohama 240-8501, Japan}
\email{nishimura-takashi-yx@ynu.jp}
\begin{document}
\maketitle
\begin{abstract}
This is a survey article on distance-squared mappings and their related topics.   
\end{abstract}
\newtheorem{theorem}{Theorem}[section] 
\newtheorem{corollary}{Corollary}[section] 
\newtheorem{lemma}{Lemma}[section]
\newtheorem{proposition}{Proposition}[section]
\newtheorem{definition}{Definition}[section]
\newtheorem{example}{Example}[section] 
\newtheorem{problem}{Problem}[section]
\newtheorem{conjecture}{Conjecture}[section]
\newcommand{\R}{\mathbb{R}}
\noindent
\section{Distance-squared mappings} 
\label{distance-squared}
Distance-squared mappings were firstly investigated in \cite{ichikinishimura}.   
\par 
Let $n$ (resp., $\mathbb{R}^n$) be a positive integer (resp., the $n$-dimensional Euclidean space). 
The $n$-dimensional Euclidean distance is the function 
$d : \mathbb{R}^n\times \mathbb{R}^n\rightarrow \mathbb{R}$ 
defined by 
\[
d(x,y)=\sqrt{\sum_{i=1}^n(x_i-y_i)^2}, 
\] 
where $x=(x_1,\ldots,x_n)$ and $y=(y_1,\ldots,y_n)$.
For any point $p\in \mathbb{R}^n$, the mapping 
$d_p : \mathbb{R}^n\rightarrow \mathbb{R}$, defined by $d_p(x)=d(p,x)$, is called a {\it distance function}.
\begin{definition}
{\rm Let $p_1,\ldots, p_\ell$ $(\ell\ge 1)$ be given points in $\mathbb{R}^n$. Then, the following mapping 
$d_{(p_1,\ldots,p_\ell)} : \mathbb{R}^n \rightarrow \mathbb{R}^\ell$ is called a {\it distance mapping}:
\[d_{(p_1,\ldots,p_\ell)}(x)=(d(p_1,x),\ldots,d(p_\ell,x)).\] }
\end{definition}
\noindent 
A distance mapping is one in which each component is a distance function.      
Distance mappings were firstly studied in the undergraduate-thesis of the first author, 
and the main result of the thesis 
is the following proposition.    
Proposotion \ref{d} can be found also in \cite{ichikinishimura} with a rigorous proof.   
The proof uses several geometric results in \cite{berger}.    
Let $S^n$ be the $n$-dimensional unit sphere in $\mathbb{R}^{n+1}$
\begin{proposition}[\cite{ichikinishimura}] \label{d}
Let $i : S^1\rightarrow i(S^1)\subset \mathbb{R}^2$ be a homeomorphism.
Then, there exist two points $p_1$,$p_2$ $\in i(S^1)$ such that $d_{(p_1,p_2)}\circ i : S^1 \rightarrow \mathbb{R}^2$ 
is homeomorphic to the image $d_{(p_1,p_2)}\circ i(S^1)$.
\end{proposition}
Proposition \ref{d} is applicable even if a mapping $i$ is not differentiable 
anywhere.    
However, it seems quite difficult to derive higher-dimensional extensions of the proposition.   
\par 
On the other hand, it is possible to obtain 
the differentiable version of higher-dimensional extensions as follows. 
\begin{definition}
{\rm Let $p_i$ $(1\le i\le \ell)$ be a given point in $\mathbb{R}^n$.   
Then, the following mapping $D_{(p_1,\ldots,p_\ell)} : \mathbb{R}^n \rightarrow \mathbb{R}^\ell$ is called a {\it distance-squared mapping}}:
\[D_{(p_1,\ldots,p_\ell)}(x)=(d^2(p_1,x),\ldots,d^2(p_\ell,x)).\]
\end{definition}
Although $D_{(p_1,\ldots,p_\ell)}$ always has a singular point, 
the following Theorems \ref{D1} and \ref{D2} hold as follows.   
\begin{theorem}[\cite{ichikinishimura}] \label{D1}
Let $M$ be an $m$-dimensional closed $C^\infty$ manifold $(m\geq 1)$, and let $i : M\to \mathbb{R}^\ell$ 
$(m+1\leq \ell)$ be a $C^\infty$ embedding. Then, there exist $p_1,\ldots,p_{m+1}\in i(M)$, $p_{m+2},\ldots,p_\ell\in \mathbb{R}^\ell $such that $D_{(p_1,\ldots,p_\ell)}\circ i : M\rightarrow \mathbb{R}^\ell$ is 
a $C^\infty$ embedding.
\end{theorem}
For the definition of embedding, see \cite{hirsh}.
\begin{corollary}[\cite{ichikinishimura}]
Let $M$ be an $m$-dimensional closed $C^\infty$ manifold $(m\geq 1)$, and let $i : M\rightarrow \mathbb{R}^{m+1}$ be a $C^\infty$ embedding. Then, there exist $p_1,\ldots,p_{m+1}\in i(M)$ 
such that $D_{(p_1,\ldots,p_{m+1})}\circ i : M\rightarrow \mathbb{R}^{m+1}$ is a $C^\infty$ embedding.
\end{corollary}
Let $M$ be an $m$-dimensional closed $C^\infty$ manifold $(m\geq 1)$, and let $i : M\rightarrow \mathbb{R}^{m+1}$ be a $C^\infty$ embedding. Then, there exist $p_1,\ldots,p_{m+1}\in i(M)$ 
such that $D_{(p_1,\ldots,p_{m+1})}\circ i : M\rightarrow \mathbb{R}^{m+1}$ is a $C^\infty$ embedding.
Let $M$ and $N$ be $C^\infty$ manifolds.   
A $C^\infty$ immersion $f : M\rightarrow N$ (for the definition of immersion, see \cite{hirsh}) is said to be 
{\it with normal crossing at a point $y\in N$} if $f^{-1}(y)$ is a finite set $\{x_1,x_2,\ldots,x_n\}$ and for any subsets $\{\lambda _1,\lambda _2,\ldots,\lambda _s\}\subset 
\{1,2,\ldots,n\}$ $(s\leq n)$, 
\[{\rm codim}\,\Bigl(\bigcap _{j=1}^sdf_{x_{\lambda _{j}}}(T_{x_{\lambda _{j}}}M)\Bigr)=\sum_{j=1}^s{\rm codim}\,(df_{x_{\lambda _{j}}}(T_{x_{\lambda _{j}}}M)),\]
where ${\rm codim}\,H={\rm dim}\,T_yN-{\rm dim}\,H$ for a linear subspace $H\subset\,T_yN$. 
A $C^\infty$ immersion $f : M\rightarrow N$ is said to be {\it with normal crossing} if $f$ is a $C^\infty$ immersion with normal crossing at any point $y\in N$.
\begin{theorem}[\cite{ichikinishimura}] \label{D2}
Let $M$ be an $m$-dimensional closed $C^\infty$ manifold $(m\geq 1)$, and let $i : M\rightarrow \mathbb{R}^\ell$ $(m+1\leq \ell)$ be a $C^\infty$ immersion with normal crossing. Then, there exist $p_1,\ldots,p_{m+1}\in i(M)$, $p_{m+2},\ldots,p_\ell\in \mathbb{R}^\ell$
such that $D_{(p_1,\ldots,p_\ell)}\circ i : M\rightarrow \mathbb{R}^\ell$ is a $C^\infty$ immersion with normal crossing.
\end{theorem} 
\begin{corollary}[\cite{ichikinishimura}]
Let $M$ be an $m$-dimensional closed $C^\infty$ manifold $(m\geq 1)$, and let 
$i : M\rightarrow \mathbb{R}^{m+1}$ be a $C^\infty$ immersion with normal crossing. Then, there exist $p_1,\ldots,p_{m+1}\in i(M)$ such that $D_{(p_1,\ldots,p_{m+1})}\circ i : M\rightarrow \mathbb{R}^{m+1}$ is 
a $C^\infty$ immersion with normal crossing.
\end{corollary}
\par 
We say that $\ell$-points $p_1,\ldots,p_\ell\in  \mathbb{R}^n$ $(1\leq \ell\leq n+1)$ are 
{\it in general position}  if $\ell=1$ or 
$\overrightarrow{p_1p_2},\ldots,\overrightarrow{p_1p_\ell}$ $(2\leq \ell\leq n+1)$ are linearly independent, 
where $\overrightarrow{p_ip_j}$ stands for $(p_{j1}-p_{i1},\ldots,p_{jn}-p_{in})$ 
$(p_i=(p_{i1},\ldots,p_{in}), p_j=(p_{j1},\ldots,p_{jn})\in \mathbb{R}^n)$.
A mapping $f : \R^n\rightarrow \mathbb{R}^\ell$ is said to be {\it $\mathcal{A}$-equivalent} to a mapping $g : \R^n\rightarrow \R^\ell$ 
if there exist $C^\infty$ diffeomorphisms $\varphi : \mathbb{R}^n\rightarrow \R^n$ and $\psi : \R^\ell\rightarrow \R^\ell$ 
such that $\psi \circ f\circ \varphi =g$.
For any two positive integers $\ell, n$ satisfying $\ell \le n$, the following mapping 
$\Phi_\ell : \R^n\rightarrow \R^\ell$ $(\ell\leq n)$ is called 
{\it the normal form of definite fold mappings}: 
\[
\Phi_\ell(x_1,\ldots,x_n)=(x_1,\ldots,x_{\ell-1}, x^2_\ell+\cdots +x^2_n).
\]
The properties of distance-squared mappings, especially (I) and (II) of the following Theorem 
\ref{main1}, are essential in the proofs of Theorems 
\ref{D1} and \ref{D2}.     
Thus, in this sense, the following Theorem \ref{main1} may be regarded as the main result in \cite{ichikinishimura}.   
\begin{theorem}[\cite{ichikinishimura}]\label{main1}
\ \\
$(I)$\ Let $\ell$,$n$ be integers such that $2\leq \ell\leq n$, and let $p_1,\ldots,p_\ell\in \mathbb{R}^n$ be in general position. Then, $D_{(p_1,\ldots,p_\ell)} : \R^n\rightarrow \R^\ell$ is $\mathcal{A}$-equivalent 
to the normal form of definite fold mappings.\\
$(II)$\ Let $\ell$,$n$ be integers such that $2\leq n<\ell $, and let $p_1,\ldots,p_{n+1}\in \mathbb{R}^n$ be in general position. Then, $D_{(p_1,\ldots,p_\ell)} : \R^n\rightarrow \R^\ell$ is $\mathcal{A}$-equivalent to the inclusion $(x_1,\ldots,x_n)\mapsto (x_1,\ldots,x_n,0,\ldots,0)$.
\end{theorem} 
All results in this section have been rigorously proved in \cite{ichikinishimura}.      
\section{Lorentzian distance-squared mappings}\label{Lorentzian}
Lorentzian distance-squared mappings were firstly studied in \cite{ichikinishimura2}.   
\par 
As same as in Section \ref{distance-squared}, we let $n$ be a positive integer.   
Let $x, y$ be two vectors of $\mathbb{R}^{n+1}$.       
Then, the {\it Lorentzian inner product} is the following qudratic form: 
\[ \langle x, y\rangle =-x_0y_0+x_1 y_1+\cdots +x_n y_n,\]
where $x=(x_0, x_1, \ldots, x_n), y=(y_0, y_1, \ldots, y_n)$.    
If the role of the Euclidean inner product $x\cdot y=\sum_{i=0}^nx_iy_i$ is replaced by 
the Lorentzian inner product, then the $(n+1)$-dimensional vector space $\mathbb{R}^{n+1}$ is called  
{\it Lorentzian $(n+1)$-space}, and it is denoted by $\mathbb{R}^{1,n}$.    
For a vector $x$ of 
Lorentzian $(n+1)$-space $\mathbb{R}^{1,n}$, {\it Lorentzian length} of $x$ is $\sqrt{\langle x,x\rangle}$.   
Notice that a pure imaginary value may be taken as the Lorentzian length 
and thus $\sqrt{\langle x,x\rangle}$ does not give a real-valued function.   
On the other hand, 
its square $x\mapsto \langle x, x\rangle$ is always a real value.  
For a non-zero vector $x\in \mathbb{R}^{1,n}$, it is called  
{\it space-like}, {\it light-like} or {\it time-like} if 
its Lorentzian length 
is positive, zero or pure imaginary respectively. 
The following is the definition of the likeness of the vector subspace.       
\begin{definition}[\cite{ratcliffe}]\label{likeness}
{\rm 
Let $V$ be a vector subspace of $\mathbb{R}^{1,n}$.   
Then $V$ is said to be 
\begin{enumerate}
\item {\it time-like} if $V$ has a time-like vector, 
\item {\it space-like} if every nonzero vector in $V$ is space-like, or 
\item {\it light-like} otherwise.    
\end{enumerate}
}
\end{definition}
\noindent 
The {\it light cone} of Lorentzian $(n+1)$-space $\mathbb{R}^{1,n}$, 
denoted by $LC$, is the set of  $x\in \mathbb{R}^{1,n}$ such that 
$\langle x, x\rangle =0$.   
For more details on Lorentzian space, 
refer to \cite{ratcliffe}.    
Recently, Singularity Theory has been very actively applied to 
geometry of submanifolds in Lorentzian space 
(for instance, see 
\cite{izumiyakossowskipeicarmen,izumiyanunocarmen,izumiyapeisano,
izumiyasaji,izumiyasato,izumiyatakahashitari,izumiyatari1,izumiyatari2,
izumiyahandan,kasedou,sato,tari1,tari2}).     
In \cite{ichikinishimura2}, 
it is given a different application of Singularity Theory to the study of Lorentzian space 
from these researches.
\par 
Let $p$ be a point of $\mathbb{R}^{1,n}$.    
The {\it Lorentzian distance-squared function} is 
the following function  
$\ell_p^2: \mathbb{R}^{1,n}\to \mathbb{R}$ (\cite{izumiyakossowskipeicarmen}):   
\[
\ell_p^2(x)=\langle x-p, x-p\rangle .   
\]
Let $p_0, \ldots, p_{k}\in \mathbb{R}^{1,n}$ $(1\le k)$ be finitely many points.    
For any $p_0, \ldots, p_{k}\in \mathbb{R}^{1,n}$, the {\it Lorentzian distance-squared mapping}, denoted by  
$L_{(p_0, \ldots, p_k)}: \mathbb{R}^{1,n}\to \mathbb{R}^{k+1}$, is defined as follows:   
\[
L_{(p_0, \ldots, p_k)}(x)=\left(\ell_{p_0}^2(x), \ldots, \ell_{p_k}^2(x)\right).   
\]
For finitely many points $p_0, \ldots, p_{k}\in \mathbb{R}^{1,n}$ $(1\le k)$, 
a vector subspace $V$ is called the {\it recognition subspace} and is denoted by 
$V(p_0, \ldots, p_k)$ of $\mathbb{R}^{1,n}$ if  the following is satisfied:  
\[
V= 
\sum_{i=1}^k \mathbb{R}\;\overrightarrow{p_0 p_i}.    
\]

For any two positive integers $k, n$ satisfying $k < n$, 
the following mapping 
$\Psi_k: \mathbb{R}^{1,n}\to \mathbb{R}^{k+1}$ is called 
the {\it normal form of Lorentzian indefinite fold mapping}: 
\[
\Psi_k\left(x_0, x_1, \ldots, x_n\right)=
\left(x_1, \ldots, x_k, 
-x_0^2+\sum_{i=k+1}^{n}x_i^2\right).
\]
Let $j, k$ be two positive integers satisfying $j\leq k$ and 
let $\tau _{(j,k)}: \mathbb{R}^{j+1}\rightarrow \mathbb{R}^{k+1}$ be the 
inclusion:
\[\tau _{(j,k)}(X_0,X_1,\ldots,X_j)=(X_0,X_1,\ldots,X_j,0,\ldots,0).\]
\begin{theorem}[\cite{ichikinishimura2}]\label{not in general position}
\begin{enumerate}
\item \label{1}
Let $k,n$ be two positive integers and let $p_0,\ldots,p_k\in \mathbb{R}^{n,1}$ 
be the same point $($i.e. $\dim V(p_0, \ldots, p_k)=0)$. Then, the Lorentzian distance-squared mapping 
$L_{(p_0,\ldots,p_k)}:\mathbb{R}^{n,1}\rightarrow \mathbb{R}^{k+1}$ is 
$\mathcal{A}$-equivalent to the mapping
\[(x_0,\ldots,x_n)\mapsto \left(-x_0^2+\sum_{i=1}^nx_i^2,0,\ldots,0\right).\]
\item \label{2}
Let $j, k, n$ be three positive integers satisfying $j<n, j\leq k$, 
and let $p_0, \ldots, p_k$ $\in \mathbb{R}^{1,n}$ be $(k+1)$ 
points such that two recognition subspaces 
$V(p_0, \ldots, p_k)$ and $V(p_0, \ldots, p_j)$ 
have the same dimension $j$. 
Then, the following hold:
\begin{enumerate}
\item \label{2.1}
The mapping 
$L_{(p_0, \ldots, p_k)}: \mathbb{R}^{1,n}\to \mathbb{R}^{k+1}$ is 
$\mathcal{A}$-equivalent to $\tau _{(j,k)}\circ \Phi_{j+1}$ if and only if 
$V(p_0, \ldots, p_k)$ is time-like. 
\item \label{2.2}
The mapping 
$L_{(p_0, \ldots, p_k)}:\mathbb{R}^{1,n}\to \mathbb{R}^{k+1}$ is 
$\mathcal{A}$-equivalent to $\tau _{(j,k)}\circ \Psi_j$ if and only if 
$V(p_0, \ldots, p_k)$ is space-like. 
\item \label{2.3}
The mapping 
$L_{(p_0, \ldots, p_k)}:\mathbb{R}^{1,n}\to \mathbb{R}^{k+1}$ is 
$\mathcal{A}$-equivalent to
\[
(x_0, \ldots, x_n)\mapsto \left( x_1, \ldots, x_j, x_0x_1
+\sum_{i=j+1}^n x_i^2,0,\ldots,0\right)
\]
if and only if 
$V(p_0, \ldots, p_k)$ is light-like. 
\end{enumerate}
\item \label{3}
Let $k, n$ be two positive integers satisfying $n\leq k$ and 
let $p_0, \ldots, p_k$ $\in \mathbb{R}^{1,n}$ be $(k+1)$ 
points such that $\dim V(p_0,\ldots,p_k)=\dim V(p_0,\ldots,p_n)=n$.
Then, the following hold: 
\begin{enumerate}
\item \label{3.1}
The mapping 
$L_{(p_0, \ldots, p_k)}:\mathbb{R}^{1,n}\to \mathbb{R}^{k+1}$ is 
$\mathcal{A}$-equivalent to $\tau _{(n,k)}\circ \Phi_{n+1}$ if and only if 
$V(p_0, \ldots, p_k)$ is time-like or space-like. 
\item \label{3.2}
The mapping 
$L_{(p_0, \ldots, p_k)}: \mathbb{R}^{1,n}\to \mathbb{R}^{k+1}$ is 
$\mathcal{A}$-equivalent to 
\[(x_0, \ldots, x_n)\mapsto (x_1, \ldots, x_n, x_0x_1,0\ldots,0)\]
if and only if 
$V(p_0, \ldots, p_k)$ is light-like.
\end{enumerate}
\item \label{4}
Let $k, n$ be two positive integers satisfying $n<k$ and 
let $p_0, \ldots, p_{k}\in \mathbb{R}^{1,n}$ be $(k+1)$ points such that 
$\dim V(p_0,\ldots,p_k)=\dim V(p_0,\ldots,p_{n+1})$
\\$=n+1$. 
Then, $L_{(p_0, \ldots, p_k)}: \mathbb{R}^{1,n}\to \mathbb{R}^{k+1}$ is always 
$\mathcal{A}$-equivalent to the inclusion 
$(x_0, \ldots, x_n)\mapsto (x_0, \ldots, x_n, 0, \ldots, 0)$. 
\end{enumerate}
\end{theorem}
Let $p_0, \ldots, p_k\in \mathbb{R}^{n,1}$ be given $(k+1)$ points.   
We say that $p_0, \ldots, p_k$ are  
{\it in general position} if 
the dimension of $V(p_0, \ldots, p_k)$ is $k$.   
For $(k+1)$ points $q_0, \ldots, q_k\in \mathbb{R}^{1,n}$ in general position $(k\le n)$, 
the singular set of $L_{(q_0, \ldots, q_k)}: \mathbb{R}^{1,n}\to \mathbb{R}^{k+1}$ is clearly the $k$-dimensional 
affine subspace spanned by these points. 
Since $\tau_{(k.k)}$ is the identity mapping, we have the following corollary.
\begin{corollary}[\cite{ichikinishimura2}]\label{corollary 2.1}
\begin{enumerate}
\label{C1} \item
Let $k, n$ be two positive integers satisfying $k< n$ 
and let $p_0, \ldots, p_k$ belonging to 
$\mathbb{R}^{1,n}$ be $(k+1)$ points in general position. 
Then, the following hold:   
\begin{enumerate}
\item The mapping 
$L_{(p_0, \ldots, p_k)}: \mathbb{R}^{1,n}\to \mathbb{R}^{k+1}$ is 
$\mathcal{A}$-equivalent to $\Phi_{k+1}$ if and only if $V(p_0, \ldots, p_k)$ is time-like. 
\item The mapping 
$L_{(p_0, \ldots, p_k)}:\mathbb{R}^{1,n}\to \mathbb{R}^{k+1}$ is 
$\mathcal{A}$-equivalent to $\Psi_k$ if and only if 
$V(p_0, \ldots, p_k)$ is space-like. 
\item The mapping 
$L_{(p_0, \ldots, p_k)}:\mathbb{R}^{1,n}\to \mathbb{R}^{k+1}$ is 
$\mathcal{A}$-equivalent to
\[
(x_0, \ldots, x_n)\mapsto \left( x_1, \ldots, x_k, x_0x_1
+\sum_{i=k+1}^n x_i^2\right)
\]
if and only if 
$V(p_0, \ldots, p_k)$ is light-like. 
\end{enumerate}
\item 
Let $n$ be a positive integer and 
let $p_0, \ldots, p_n\in \mathbb{R}^{1,n}$ be $(n+1)$ points in general position. 
Then, the following hold:   
\begin{enumerate}
\item The mapping 
$L_{(p_0, \ldots, p_n)}:\mathbb{R}^{1,n}\to \mathbb{R}^{n+1}$ is 
$\mathcal{A}$-equivalent to $\Phi_{n+1}$ if and only if 
$V(p_0, \ldots, p_n)$ is time-like or space-like. 
\item The mapping 
$L_{(p_0, \ldots, p_n)}: \mathbb{R}^{1,n}\to \mathbb{R}^{n+1}$ is 
$\mathcal{A}$-equivalent to
\[(x_0, \ldots, x_n)\mapsto (x_1, \ldots, x_n, x_0x_1)\]
if and only if 
$V(p_0, \ldots, p_n)$ is light-like.   
\end{enumerate}    
\end{enumerate}
\end{corollary}

The following are clear:   
\begin{enumerate} 
\item Any non-singular fiber of $\Phi_{n}$ is a circle.     
\item Any non-singular fiber of $\Psi_{n-1}$ is an equilateral hyperbola.   
\item Any non-singular fiber of $(x_0, \ldots, x_n)\mapsto (x_1, \ldots, x_{n-1}, x_0x_1+x_n^2)$ is 
a parabola (possibly at infinity).    
\end{enumerate}
Therefore, by the case $k=n-1$ in Corollary \ref{C1}, we have the following:
\begin{corollary}[\cite{ichikinishimura2}]\label{corollary 1}
Let $n$ be a positive integer such that $\; 2\le n$ and 
let $p_0, \ldots, p_{n-1}$ belonging to 
$\mathbb{R}^{1,n}$ be $n$ points in general position. 
Then, the following hold:   
\begin{enumerate}
\item There exists a $C^\infty$ diffeomorphism $h:\mathbb{R}^{1,n}\to \mathbb{R}^{1,n}$ 
by which any non-singular fiber $L_{(p_0, \ldots, p_{n-1})}^{-1}(y)$ is mapped to a circle 
if and only if the recognition subspace $V(p_0, \ldots, p_{n-1})$ is time-like.   
\item There exists a $C^\infty$ diffeomorphism $h:\mathbb{R}^{1,n}\to \mathbb{R}^{1,n}$ 
by which any non-singular fiber $L_{(p_0, \ldots, p_{n-1})}^{-1}(y)$ is mapped to 
an equilateral hyperbola 
if and only if $V(p_0, \ldots, p_{n-1})$ is space-like.  
\item There exists a $C^\infty$ diffeomorphism $h:\mathbb{R}^{1,n}\to \mathbb{R}^{1,n}$ 
by which any non-singular fiber $L_{(p_0, \ldots, p_{n-1})}^{-1}(y)$ is mapped to a parabola 
if and only if the recognition subspace  
$V(p_0, \ldots, p_{n-1})$ is light-like.  
\end{enumerate}
\end{corollary}
\noindent 
In \cite{ichikinishimura2}, it is remarked that an affine diffeomorphism can be chosen 
as the diffeomorphism $h:\mathbb{R}^{1,n}\to \mathbb{R}^{1,n}$ in Corollary \ref{corollary 1}.    
\par 
\smallskip 
The motivation to classify Lorentzian distance-squared mappings in \cite{ichikinishimura2} 
is the classification results on 
distance-squared mappings, namely Theorem \ref{main1}.  
It is natural to ask how Theorem \ref{main1} changes 
if distance-squared functions are replaced with Lorentzian 
distance-squared functions.
\noindent 
Combining Theorem \ref{main1} and Corollary \ref{corollary 2.1}, we have  
the following:  
\begin{corollary}[\cite{ichikinishimura2}]\label{corollary 2}
\begin{enumerate}
\item 
Let $k, n$ be two positive integers satisfying $ k< n$ and 
let $p_0, \ldots, p_k$ belonging to $\mathbb{R}^{1,n}$ be 
$(k+1)$ points in general position. 
Then, $L_{(p_0, \ldots, p_k)}$ 
is $\mathcal{A}$-equivalent to 
$D_{(p_0, \ldots, p_k)}$ 
if and only if 
$V(p_0, \ldots, p_k)$ is time-like. 
\item 
Let $n$ be a positive integer and 
let $p_0, \ldots, p_n\in \mathbb{R}^{1,n}$ be $(n+1)$ points in general position. 
Then, $L_{(p_0, \ldots, p_n)}$ 
is $\mathcal{A}$-equivalent to 
$D_{(p_0, \ldots, p_n)}$ 
if and only if 
$V(p_0, \ldots, p_n)$ is time-like or space-like. 
\item 
Let $k, n$ be two positive integers satisfying $n<k$ and 
let $p_0, \ldots, p_{k}\in \mathbb{R}^{1,n}$ be $(k+1)$ points such that 
the $(n+2)$ points $p_0, \ldots, p_{n+1}$ are in general position.    
Then, $L_{(p_0, \ldots, p_k)}$ 
is always $\mathcal{A}$-equivalent to 
$D_{(p_0, \ldots, p_k)}$.   
\end{enumerate}
\end{corollary}
All results in this section have been rigorously proved in \cite{ichikinishimura2}.   
\section{Generalized distance-squared mappings of $\mathbb{R}^2$ into $\mathbb{R}^2$}\label{section 3}
Generalized distance-squared mappings were firstly studied in \cite{ichikinishimura3}.   
For any two positive integers $k, n$, we let
$p_0, p_1, \ldots, p_k$ be $(k+1)$ points of $\mathbb{R}^{n+1}$.
We set $p_i=(p_{i0}, p_{i1}, \ldots, p_{in})$ $(0\le i\le k)$.
We let $A=(a_{ij})_{0\le i\le k, 0\le j\le n}$ be a $(k+1)\times (n+1)$ matrix with non-zero entries.
Then,  we consider the following mapping $G_{(p_0, p_1, \ldots, p_k, A)}: \mathbb{R}^{n+1}\to \mathbb{R}^{k+1}$: 
{\small 
\[
G_{(p_0, p_1, \ldots, p_k, A)}(x)=\left(
\sum_{j=0}^n a_{0j}(x_j-p_{0j})^2,
\sum_{j=0}^n a_{1j}(x_j-p_{1j})^2,
\ldots,
\sum_{j=0}^n a_{{k}j}(x_j-p_{{k}j})^2
\right),
\]
}
where $x=(x_0, x_1, \ldots, x_n)$
The mapping $G_{(p_0, p_1, \ldots, p_k, A)}$ is called
a {\it generalized distance-squared mapping}.   
Notice that a distance-squared mapping 
$D_{(p_0, p_1, \ldots, p_k)}$ defined in Section \ref{distance-squared} is 
the mapping $G_{(p_0, p_1,
\ldots, p_k,A)}$ in the case that each entry of $A$ is $1$, and 
a Lorentzian distance-squared mapping 
$L_{(p_0, p_1, \ldots, p_k)}$ defined in Section \ref{Lorentzian} is 
the mapping $G_{(p_0, p_1, \ldots, p_k,A)}$ in the case of $a_{i0}=-1 $ and $a_{ij}=1$ if $j\ne 0$.   
Notice also that in these
cases, the rank of $A$ is $1$.    
In the applications of singularity theory to
differential geometry, generalized distance-squared mappings
are a useful tool.   Information on the
contacts amongst the families of quadrics defined by the components
of $G_{(p_0,p_1,\ldots,p_k,A)}$ 
is given by their singularities.    
Hence, it is natural to classify
maps $G_{(p_0,p_1,\ldots,p_k,A)}$.   
\par 
\medskip 
From now on in this section, we concentrate on the case of $\mathbb{R}^2$ into $\mathbb{R}^2$.   
It is interesting to observe that new $\mathcal{A}$-classes occur even in this case.
\begin{definition}\label{definition 1}
{\rm
\begin{enumerate}
\item
Let $\Phi_{n+1}: \mathbb{R}^{n+1}\to \mathbb{R}^{n+1}$ denote the following mapping:
\[
\Phi_{n+1}(x_0, x_1, \ldots, x_n) =\left(x_0, x_1, \ldots, x_{n-1}, x_n^2\right).
\]
When a map-germ $f: (\mathbb{R}^{n+1}, q)\to (\mathbb{R}^{n+1}, f(q))$ is $\mathcal{A}$-equivalent to
$\Phi_{n+1}: (\mathbb{R}^{n+1}, 0)\to (\mathbb{R}^{n+1},0)$,
the point $q\in \mathbb{R}^{n+1}$ is said to be a {\it fold point of $f$}.
\item
Let $\Gamma_{n+1}: \mathbb{R}^{n+1}\to \mathbb{R}^{n+1}$ denote the following mapping:
\[
\Gamma_{n+1}(x_0, x_1, \ldots, x_n) =\left(x_0, x_1, \ldots, x_{n-1}, x_n^3+x_0x_n\right).
\]
When a map-germ $f: (\mathbb{R}^{n+1}, q)\to (\mathbb{R}^{n+1}, f(q))$ is $\mathcal{A}$-equivalent to
$\Gamma_{n+1}: (\mathbb{R}^{n+1}, 0)\to (\mathbb{R}^{n+1},0)$,
the point $q\in \mathbb{R}^{n+1}$ is said to be a {\it cusp point of $f$}.
\end{enumerate}
}
\end{definition}
\noindent
It is known that both $\Phi_{n+1}, \Gamma_{n+1}$ are proper and stable mappings (for instance see \cite{arnoldetall}).
\par
Recall the special cases of Theorem \ref{main1} and Corollary \ref{corollary 2.1} as follows:   
\begin{proposition}[special cases of Theorem \ref{main1} and Corollary \ref{corollary 2.1}]
\label{proposition 3.1}
Let $p_0, p_1, \ldots, p_n$ be $(n+1)$-points of $\mathbb{R}^{n+1}$ such that the dimension of $\sum_{i=1}^n\mathbb{R}\overrightarrow{p_0p_i}$ is $n$.
Then, the following hold:
\begin{enumerate}
\item The distance-squared mapping $D_{(p_0, p_1, \ldots, p_n)}: \mathbb{R}^{n+1}\to \mathbb{R}^{n+1}$
is $\mathcal{A}$-equivalent to $\Phi_{n+1}$.
\item The Lorentzian distance-squared mapping $L_{(p_0, p_1, \ldots, p_n)}: \mathbb{R}^{n+1}\to \mathbb{R}^{n+1}$
is $\mathcal{A}$-equivalent to $\Phi_{n+1}$.
\end{enumerate}
\end{proposition}
For generalized distance-squared mappings,
it is natural to expect that for generic $p_0, p_1, \ldots, p_n$, $G_{(p_0, p_1, \ldots, p_n, A)}$ is proper and stable, and the rank
of $A$ is a complete invariant of $\mathcal{A}$-types.   
Thus, we reach the following conjecture.  
\begin{conjecture}[\cite{ichikinishimura3}]\label{conjecture 1}
Let $A_k$
be an $(n+1)\times (n+1)$ matrix of rank $k$ with non-zero entries $(1\le k\le (n+1))$.
Then, there exists a subset $\Sigma\subset (\mathbb{R}^{n+1})^{n+1}$ of Lebesgue measure zero such that
for any $(p_0, p_1, \ldots, p_n)\in  (\mathbb{R}^{n+1})^{n+1}- \Sigma$, the following hold:
\begin{enumerate}
\item For any $k$ $(1\le k\le (n+1))$,
the generalized distance-squared mapping $G_{(p_0, p_1, \ldots, p_n, A_k)}$ is proper and stable.
\item For any two integers $k_1, k_2$ such that $1\le k_1<k_2\le (n+1)$,
$G_{(p_0, p_1, \ldots, p_n, A_{k_2})}$ is not $\mathcal{A}$-equivalent to $G_{(p_0, p_1, \ldots, p_n, A_{k_1})}$.
\item Let $B_k$ be an $(n+1)\times (n+1)$ matrix of rank $k$ with non-zero entries $(1\le k\le (n+1))$ and
let  $(q_0, q_1, \ldots, q_n)$ be in $ (\mathbb{R}^{n+1})^{n+1}- \Sigma$.     Then,
$G_{(p_0, p_1, \ldots, p_n, A_k)}$ is $\mathcal{A}$-equivalent to $G_{(q_0, q_1, \ldots, q_n, B_k)}$ for any $k$.
\end{enumerate}
\end{conjecture}
In \cite{ichikinishimura3}, the affirmative answer to Conjecture \ref{conjecture 1} in the case $n=1$ are given  as follows: 
\begin{theorem}[\cite{ichikinishimura3}] \label{theorem 1}
Let $((x_0, y_0), (x_1, y_1))$ be the standard coordinates of $(\mathbb{R}^2)^2$
and let $\Sigma$ be the hypersurface in $(\mathbb{R}^2)^2$ defined by
$(x_0-x_1)(y_0-y_1)=0$.
Let $(p_0, p_1)$ be a point in $(\mathbb{R}^2)^2 - \Sigma$ and
let $A_k$ be a $2\times 2$ matrix of rank $k$ with non-zero entries (k=1, 2).
Then, the following hold:
\begin{enumerate}
\item
The mapping $G_{(p_0, p_1, A_1)}$ is $\mathcal{A}$-equivalent to $\Phi_2$.
\item
The mapping $G_{(p_0, p_1, A_2)}$ is proper and stable, and
it is not $\mathcal{A}$-equivalent to  $G_{(p_0, p_1, A_1)}$.
\item Let $B_2$ be a $2\times 2$ matrix of rank $2$ with non-zero entries and let
$(q_0, q_1)$ be a point in $(\mathbb{R}^2)^2 - \Sigma$.
Then, $G_{(p_0, p_1, A_2)}$ is $\mathcal{A}$-equivalent
to $G_{(q_0, q_1, B_2)}$.
\end{enumerate}
\end{theorem}
\noindent
{ 
There is another motivation for Theorem \ref{theorem 1}.    
Set $f_t(x)=x+t x^2$ $(t, x\in \mathbb{R})$.   Then, the following two 
are easily observed.   
\begin{enumerate} 
\item $f_t$ is proper and stable for any $t\in \mathbb{R}$.   
\item $f_t$ $(t\ne 0)$ is not $\mathcal{A}$-equivalent to $f_0$.    
\end{enumerate}
Notice that the mapping  
\[
F: 
\mathbb{R}\to C^\infty(\mathbb{R}, \mathbb{R})
\]
defined by $F(t)=f_t$, is continuous nowhere.   Here,  
$C^\infty(\mathbb{R}, \mathbb{R})$ is the topological space consisting of 
$C^\infty$ mappings $\mathbb{R}\to \mathbb{R}$ endowed with 
the Whitney $C^\infty$ topology.       
The one-parameter family $f_t$ is a very simple example for 
preliminary phenomena of wall crossing phenomena.   
Theorem \ref{theorem 1} gives an example for such phenomena 
in the case of the plane to the plane as follows.    
Let $M(2, \mathbb{R})$ be the set consisting of $2\times 2$ matrices with real entries and 
let $P: \mathbb{R}\to (\mathbb{R}^2)^2 - \Sigma$ be a continuous mapping, where 
$\Sigma$ is the set given in Theorem \ref{theorem 1}.    
Moreover, let $A: \mathbb{R}\to M(2,\mathbb{R})$ be a continuous mapping such that 
rank$A(0)=1$ and rank$A(t)=2$ if $t\ne 0$.    
Then, Theorem \ref{theorem 1} implies the following interesting phenomenon:   
\begin{enumerate}
\item The mapping $G_{(P(s), A(t))}$ is proper and stable for any $(s, t)\in \mathbb{R}^2$.   
\item The mapping $G_{(P(s), A(t))}$ $(t\ne 0)$ is not $\mathcal{A}$-equivalent to 
$G_{(P(s), A(0))}$.   
\end{enumerate}
Notice that the mappings $P$ and $A$ induce the mapping 
\[
(\widetilde{P}, \widetilde{A}): 
\mathbb{R}^2\to C^\infty(\mathbb{R}^2, \mathbb{R}^2)
\]
defined by $(\widetilde{P}, \widetilde{A})(s, t)=G_{(P(s), A(t))}$, where 
$C^\infty(\mathbb{R}^2, \mathbb{R}^2)$ is the topological space consisting of 
$C^\infty$ mappings $\mathbb{R}^2\to \mathbb{R}^2$ endowed with 
the Whitney $C^\infty$ topology.       
Notice also that the mapping $(\widetilde{P}, \widetilde{A})$ is continuous nowhere.      
Therefore, $(\widetilde{P}, \widetilde{A})$ is useless for the proof of (3) of Theorem \ref{theorem 1}.    
}
\par 
\medskip 
The keys for proving Theorem \ref{theorem 1} are the following two propositions.
\begin{proposition}[\cite{ichikinishimura3}]\label{proposition 2}
Let $A_2$ be a $2\times 2$ matrix of rank two with non-zero entries.
Let $p_0, p_1$ be two points of $\mathbb{R}^2$
satisfying $(p_0, p_1)\in (\mathbb{R}^2)^2-\Sigma
$,
where $\Sigma\subset (\mathbb{R}^2)^2$ is the hypersurface defined in Theorem \ref{theorem 1}.
Then, the following hold:
\begin{enumerate}
\item The singular set $S(G_{(p_0, p_1, A_2)})$ is {a rectangular} hyperbola.
\item Any point of $S(G_{(p_0, p_1, A_2)})$is a fold point except for one.
\item The exceptional point given in (2) is a cusp point.
\end{enumerate}
\end{proposition}
\begin{proposition}[\cite{ichikinishimura3}]\label{proposition 3}
Let $A_2$ be a $2\times 2$ matrix of rank two with non-zero entries.
Let $p_0, p_1$ be two points of $\mathbb{R}^2$
satisfying $(p_0, p_1)\in (\mathbb{R}^2)^2-\Sigma
$,
where $\Sigma\subset (\mathbb{R}^2)^2$ is the hypersurface defined in Theorem \ref{theorem 1}.
Then, for any positive real numbers $a, b$ $(a\ne b)$,
there exists a point $q=(q_0, q_1)\in \mathbb{R}^2$ such that
$(q, (0,0))\in (\mathbb{R}^2)^2-\Sigma$ and
$G_{(p_0, p_1, A_2)}$ is $\mathcal{A}$-equivalent to $F_q: \mathbb{R}^2\to \mathbb{R}^2$ defined by
\[
F_q(x,y)=\left((x-q_0)^2+(y-q_1)^2, ax^2+by^2\right).
\]

\end{proposition}
All results in this section have been rigorously proved in \cite{ichikinishimura3}. 
\section{Generalized distance-squared mappings of $\mathbb{R}^{n+1}$ into $\mathbb{R}^{2n+1}$}\label{section 4}
Let $n$ be a positive integer.   Generalized distance-squared mappings 
$G_{(p_0, p_1, \ldots, p_{2n}, A)}: \mathbb{R}^{n+1}\to \mathbb{R}^{2n+1}$ were firstly studied in 
\cite{ichikinishimura4}.    
In this section, we survey results for generalized distance-squared mappings of $\mathbb{R}^{n+1}$ into $\mathbb{R}^{2n+1}$ obtained in \cite{ichikinishimura4}.    
 
\par 
In the case of $n=2k$,  
a partial classification result for $G_{(p_0, \ldots, p_k, A)}$ is known as follows.    
A {\it distance-squared mapping}   
$D_{(p_0, p_1, \ldots, p_k)}$ 
(resp., {\it Lorentzian distance-squared mappings}  
$L_{(p_0, p_1, \ldots, p_k)}$) 
is the mapping 
$G_{(p_0, p_1, \ldots, p_k,A)}$ in the case that each entry of $A$ is $1$ 
(resp., in the case of $a_{i0}=-1 $ and $a_{ij}=1$ if $j\ne 0$).   
\begin{proposition}[\cite{ichikinishimura, ichikinishimura2}]
\label{proposition 2}
There exists a closed subset $\Sigma\subset (\mathbb{R}^{n+1})^{2n+1}$ 
with Lebesgue measure zero 
such that for any $p=(p_0, \ldots, p_{2n})\in (\mathbb{R}^{n+1})^{2n+1}-\Sigma$,  
both $D_{(p_0, p_1, \ldots, p_{2n})}$ and $L_{(p_0, p_1, \ldots, p_{2n})}$ 
are $\mathcal{A}$-equivalent to an inclusion. 
\end{proposition}
The following Theorem \ref{theorem 4.1} generalizes Proposition \ref{proposition 2}.    
\begin{theorem}[\cite{ichikinishimura4}] 
\label{theorem 4.1}
Let $A=(a_{ij})_{0\le i\le 2n, 0\le j\le n}$ be a $(2n+1)\times (n+1)$ matrix with 
non-zero entries.   
Then, the following two hold: 
\begin{enumerate}
\item Suppose that the rank of $A$ is $n+1$.    Then, there exists a closed subset $\Sigma_A\subset (\mathbb{R}^{n+1})^{2n+1}$ 
with Lebesgue measure zero 
such that for any $p=(p_0, \ldots, p_{2n})\in  (\mathbb{R}^{n+1})^{2n+1}-\Sigma_A$, 
$G_{(p, A)}$ is $\mathcal{A}$-equivalent to the following mapping:   
\[
(x_0, x_1, \ldots, x_n)\mapsto 
(x_0^2, x_0x_1, \ldots, x_0x_n, x_1, \ldots, x_n).   
\]
\item Suppose that the rank of $A$ is less than $n+1$.    
Then, there exists a closed subset $\Sigma_A\subset (\mathbb{R}^{n+1})^{2n+1}$ with Lebesgue measure zero 
such that for any $p=(p_0, \ldots, p_{2n})\in  (\mathbb{R}^{n+1})^{2n+1}-\Sigma_A$, 
$G_{(p, A)}$ is $\mathcal{A}$-equivalent to an inclusion.      
\end{enumerate}
\end{theorem} 
The mapping given in the assertion (1) of Theorem \ref{theorem 4.1} is defined in \cite{whitneyumbrella} and is called {\it the normal form of Whitney umbrella}.     
It is clear that the normal form of Whitney umbrella is not $\mathcal{A}$-equivalent to an inclusion.    
Moreover, it is easily seen that these two mappings are proper and stable by the characterization 
theorem of stable mappings given in \cite{mather5}.         
Thus, Theorem \ref{theorem 4.1} may be regarded as a result of Theorem 
\ref{theorem 1} type.    
On the other hand, it is desirable to improve Theorem \ref{theorem 4.1} so that the bad set $\Sigma_A$ given in Theorem \ref{theorem 4.1} does not depend on the given matrix $A$.  
However, contrary to the case of $\mathbb{R}^2$ into $\mathbb{R}^2$, 
in this case it is impossible to expect the existence of such a universal bad set as follows.   
\begin{theorem}[\cite{ichikinishimura4}] 
\label{theorem 2}
There does not  exist a closed subset $\Sigma\subset (\mathbb{R}^{n+1})^{2n+1}$  with Lebesgue measure zero 
such that for any $p=(p_0, \ldots, p_{2n})\in  (\mathbb{R}^{n+1})^{2n+1}-\Sigma$
the following two hold.    
\begin{enumerate}
\item Suppose that the rank of $A$ is $n+1$.    Then, 
$G_{(p, A)}$ is $\mathcal{A}$-equivalent to the following mapping:   
\[
(x_0, x_1, \ldots, x_n)\mapsto 
(x_0^2, x_0x_1, \ldots, x_0x_n, x_1, \ldots, x_n).   
\]
\item Suppose that the rank of $A$ is less than $n+1$.    
Then, 
$G_{(p, A)}$ is $\mathcal{A}$-equivalent to an inclusion.      
\end{enumerate}
\end{theorem} 
All results in this section have been rigorously proved in \cite{ichikinishimura4}.


\end{document}